# Uniqueness and multiplicity of infinite clusters

**Geoffrey Grimmett**[1]

*University of Cambridge*

**Abstract:** The Burton–Keane theorem for the almost-sure uniqueness of infinite clusters is a landmark of stochastic geometry. Let $\mu$ be a translation-invariant probability measure with the finite-energy property on the edge-set of a $d$-dimensional lattice. The theorem states that the number $I$ of infinite components satisfies $\mu(I \in \{0,1\}) = 1$. The proof is an elegant and minimalist combination of zero–one arguments in the presence of amenability. The method may be extended (not without difficulty) to other problems including rigidity and entanglement percolation, as well as to the Gibbs theory of random-cluster measures, and to the central limit theorem for random walks in random reflecting labyrinths. It is a key assumption on the underlying graph that the boundary/volume ratio tends to zero for large boxes, and the picture for non-amenable graphs is quite different.

## 1. Introduction

The Burton–Keane proof of the uniqueness of infinite clusters is a landmark in percolation theory and stochastic geometry. The general issue is as follows. Let $\omega$ be a random subset of $\mathbb{Z}^d$ with law $\mu$, and let $I = I(\omega)$ be the number of unbounded components of $\omega$. Under what reasonable conditions on $\mu$ is it the case that: either $\mu(I = 0) = 1$, or $\mu(I = 1) = 1$? This question arose first in percolation theory with $\mu = \mu_p$, where $\mu_p$ denotes product measure (on either the vertex-set or the edge-set of $\mathbb{Z}^d$) with density $p$. It was proved in [2] that $\mu_p(I = 1) = 1$ for any value of $p$ for which $\mu_p(I \geq 1) > 0$, and this proof was simplified in [9]. Each of these two proofs utilized a combination of geometrical arguments together with a large-deviation estimate.

The true structure of the problem emerged only in the paper of Robert Burton and Michael Keane [5]. Their method is elegant and beautiful, and rests on the assumptions that the underlying measure $\mu$ is translation-invariant with a certain 'finite-energy property', and that the underlying graph is amenable (in that the boundary/volume ratio tends to zero in the limit for large boxes). The Burton–Keane method is canonical of its type, and is the first port of call in any situation where such a uniqueness result is needed. It has found applications in several areas beyond connectivity percolation, and the purpose of this paper is to summarize the method, and to indicate some connections to other problems in the theory of disordered media.

Michael Keane's contributions to the issue of uniqueness are not confined to [5]. The results of that paper are extended in [10] to long-range models (see also [33]), and to models on half-spaces. In a further paper, [6], he explored the geometrical properties of infinite clusters in two dimensions, and in [11] the existence of circuits. He wrote in the earlier paper [21] of uniqueness in long-range percolation.

[1]Statistical Laboratory, Centre for Mathematical Sciences, University of Cambridge, Wilberforce Road, Cambridge CB3 0WB, United Kingdom, e-mail: grg@statslab.cam.ac.uk

The Burton–Keane approach to uniqueness is sketched in Section 2 in the context of percolation. Its applications to rigidity percolation and to entanglement percolation are summarized in Sections 3 and 4. An application to the random-cluster model is described in Section 5, and another to random walks in random reflecting labyrinths in Section 6. The reader is reminded in Section 7 that infinite clusters may be far from unique when the underlying graph is non-amenable. We shall make periodic references to lattices, but no formal definition is given here.

## 2. Uniqueness of infinite percolation clusters

The Burton–Keane argument is easiest described in the context of percolation, and we begin therefore with a description of the bond percolation model. Let $G = (V, E)$ be a countably infinite connected graph with finite vertex-degrees. The configuration space of the model is the set $\Omega = \{0, 1\}^E$ of all 0/1-vectors $\omega = (\omega(e) : e \in E)$. An edge $e$ is called *open* (respectively, *closed*) in the configuration $\omega$ if $\omega(e) = 1$ (respectively, $\omega(e) = 0$). The product space $\Omega$ is endowed with the $\sigma$-field $\mathcal{F}$ generated by the finite-dimensional cylinder sets. For $p \in [0, 1]$, we write $\mu_p$ for product measure with density $p$ on $(\Omega, \mathcal{F})$.

The percolation model is central to the study of disordered geometrical systems, and a reasonably full account may be found in [16].

Let $\omega \in \Omega$, write $\eta(\omega) = \{e \in E : \omega(e) = 1\}$ for the set of open edges of $\omega$, and consider the open subgraph $G_\omega = (V, \eta(\omega))$ of $G$. For $x, y \in V$, we write $x \leftrightarrow y$ if $x$ and $y$ lie in the same component of $G_\omega$. We write $x \leftrightarrow \infty$ if the component of $G_\omega$ containing $x$ is infinite, and we let $\theta_x(p) = \mu_p(x \leftrightarrow \infty)$. The number of infinite components of $G_\omega$ is denoted by $I = I(\omega)$.

It is standard that, for any given $p \in [0, 1]$,

$$\text{for all } x, y \in V, \qquad \theta_x(p) = 0 \text{ if and only if } \theta_y(p) = 0, \tag{1}$$

and that

$$\theta_x(p) \begin{cases} = 0 & \text{if } p < p_c(G), \\ > 0 & \text{if } p > p_c(G), \end{cases} \tag{2}$$

where the *critical probability* $p_c(G)$ is given by

$$p_c(G) = \sup\{p : \mu_p(I = 0) = 1\}. \tag{3}$$

The event $\{I \geq 1\}$ is independent of the states of any finite collection of edges. Since the underlying measure is product measure, it follows by the Kolmogorov zero–one law that $\mu_p(I \geq 1) \in \{0, 1\}$, and hence

$$\mu_p(I \geq 1) \begin{cases} = 0 & \text{if } p < p_c(G), \\ = 1 & \text{if } p > p_c(G). \end{cases}$$

It is a famous open problem to determine for which graphs it is the case that $\mu_{p_c}(I \geq 1) = 0$, see Chapters 8–10 of [16].

We concentrate here on the case when $G$ is the $d$-dimensional hypercubic lattice. Let $\mathbb{Z} = \{\ldots, -1, 0, -1, \ldots\}$ be the integers, and $\mathbb{Z}^d$ the set of all $d$-vectors $x = (x_1, x_2, \ldots, x_d)$ of integers. We turn $\mathbb{Z}^d$ into a graph by placing an edge between any two vertices $x$, $y$ with $|x - y| = 1$, where

$$|z| = \sum_{i=1}^{d} |z_i|, \qquad z \in \mathbb{Z}^d.$$

We write $\mathbb{E}$ for the set of such edges, and $\mathbb{L}^d = (\mathbb{Z}^d, \mathbb{E})$ for the ensuing graph. Henceforth, we let $d \geq 2$ and we consider bond percolation on the graph $\mathbb{L}^d$. Similar results are valid for any lattice in two or more dimensions, and for site percolation. A *box* $\Lambda$ is a subset of $\mathbb{Z}^d$ of the form $\prod_{i=1}^{d}[x_i, y_i]$ for some $x, y \in \mathbb{Z}^d$. The *boundary* $\partial S$ of the set $S$ of vertices is the set of all vertices in $S$ which are incident to some vertex not in $S$.

A great deal of progress was made on percolation during the 1980s. Considerable effort was spent on understanding the subcritical phase (when $p < p_c$) and the supercritical phase (when $p > p_c$). It was a key discovery that, for any $p$ with $\mu_p(I \geq 1) = 1$, we have that $\mu_p(I = 1) = 1$; that is, the infinite cluster is (almost surely) unique whenever it exists.

**Theorem 1 ([2]).** *For any $p \in [0, 1]$, either $\mu_p(I = 0) = 1$ or $\mu_p(I = 1) = 1$.*

This was first proved in [2], and with an improved proof in [9]. The definitive proof is that of Burton and Keane, [5], and we sketch this later in this section. Examination of the proof reveals that it relies on two properties of the product measure $\mu_p$, namely translation-invariance and finite-energy. The first of these is standard, the second is as follows. A probability measure $\mu$ on $(\Omega, \mathcal{F})$ is said to have the *finite-energy property* if, for all $e \in \mathbb{E}$,

$$0 < \mu(e \text{ is open} \mid \mathcal{T}_e) < 1 \qquad \mu\text{-almost-surely},$$

where $\mathcal{T}_e$ denotes the $\sigma$-field generated by the states of edges other than $e$. The following generalization of Theorem 1 may be found in [5].

**Theorem 2 ([5]).** *Let $\mu$ be a translation-invariant probability measure on $(\Omega, \mathcal{F})$ with the finite-energy property. Then $\mu(I \in \{0, 1\}) = 1$.*

If, in addition, $\mu$ is ergodic, then $I$ is $\mu$-almost-surely constant, and hence: either $\mu(I = 0) = 1$ or $\mu(I = 1) = 1$. A minor complication arises for translation-invariant non-ergodic measures, and this is clarified in [6] and [12], page 42.

*Proof of Theorem 1.* The claim is trivial if $p = 0, 1$, and we assume henceforth that $0 < p < 1$. There are three steps. Since $I$ is a translation-invariant function and $\mu_p$ is ergodic, $I$ is $\mu_p$-almost-surely constant. That is, there exists $i_p \in \{1, 2, \ldots\} \cup \{\infty\}$ such that

$$\mu_p(I = i_p) = 1. \qquad (4)$$

Secondly, let us assume that $2 \leq i_p < \infty$. There exists a box $\Lambda$ such that

$$\mu_p(\Lambda \text{ intersects } i_p \text{ infinite clusters}) > 0.$$

By replacing the state of every edge in $\Lambda$ by 1, we deduce by finite-energy that

$$\mu_p(I = 1) > 0,$$

in contradiction of (4). Therefore, $i_p \in \{0, 1, \infty\}$.

In the third step we prove that $i_p \neq \infty$. Suppose on the contrary that $i_p = \infty$. We will derive a contradiction by a geometrical argument. A vertex $x$ is called a *trifurcation* if:

  (i) $x$ lies in an infinite open cluster,
  (ii) there exist exactly three open edges incident to $x$, and
  (iii) the deletion of $x$ and its three incident open edges splits this infinite cluster into exactly three disjoint infinite clusters and no finite clusters.

We write $T_x$ for the event that $x$ is a trifurcation.

By translation-invariance, the probability of $T_x$ does not depend on the choice of $x$, and thus we set $\tau = \mu_p(T_x)$. Since $i_p = \infty$ by assumption, there exists a box $\Lambda$ such that

$$\mu_p(\Lambda \text{ intersects three or more infinite clusters}) > 0.$$

On this event, we may alter the configuration inside $\Lambda$ in order to obtain the event $T_0$. We deduce by the finite-energy property of $\mu_p$ that $\tau > 0$.

The mean number of trifurcations inside $\Lambda$ is $\tau|\Lambda|$. This implies a contradiction, as indicated by the following rough argument. Select a trifurcation ($t_1$, say) of $\Lambda$, and choose some vertex $y_1$ ($\in \partial\Lambda$) which satisfies $t_1 \leftrightarrow y_1$ in $\Lambda$. We now select a new trifurcation $t_2 \in \Lambda$. By the definition of the term 'trifurcation', there exists $y_2 \in \partial\Lambda$ such that $y_1 \neq y_2$ and $t_2 \leftrightarrow y_2$ in $\Lambda$. We continue similarly, at each stage picking a new trifurcation $t_k \in \Lambda$ and a new vertex $y_k \in \partial\Lambda$. If there exist $N$ trifurcations in $\Lambda$, then we obtain $N$ distinct vertices $y_k$ lying in $\partial\Lambda$. Therefore $|\partial\Lambda| \geq N$. We take expectations to find that $|\partial\Lambda| \geq \tau|\Lambda|$, which is impossible with $\tau > 0$ for large $\Lambda$. We deduce by this contradiction that $i_p \neq \infty$. The necessary rigour may be found in [5, 16]. □

## 3. Rigidity percolation

Theorems 1 and 2 assert the almost-sure uniqueness of the infinite *connected* component. In certain other physical situations, one is interested in topological properties of subgraphs of $\mathbb{L}^d$ other than connectivity, of which two such properties are 'rigidity' and 'entanglement'. The first of these properties may be formulated as follows.

Let $G = (V, E)$ be a finite graph and let $d \geq 2$. An *embedding* of $G$ into $\mathbb{R}^d$ is an injection $f : V \to \mathbb{R}^d$. A *framework* $(G, f)$ is a graph $G$ together with an embedding $f$. A *motion* of a framework $(G, f)$ is a differentiable family $\mathbf{f} = (f_t : 0 \leq t \leq 1)$ of embeddings of $G$, containing $f$, which preserves all edge-lengths. That is to say, we require that $f = f_T$ for some $T$, and that

$$\|f_t(u) - f_t(v)\| = \|f_0(u) - f_0(v)\| \tag{5}$$

for all edges $\langle u, v \rangle \in E$, where $\|\cdot\|$ is the Euclidean norm on $\mathbb{R}^d$. We call the motion $\mathbf{f}$ *rigid* if (5) holds for *all* pairs $u, v \in V$ rather than adjacent pairs only. A framework is called *rigid* if all its motions are rigid motions.

The above definition depends on the value of $d$ and on the initial embedding $f$. For given $d$, the property of rigidity is 'generic' with respect to $f$, in the sense that there exists a natural measure $\pi$ (generated from Lebesgue measure) on the set of embeddings of $G$ such that: either $(G, f)$ is rigid for $\pi$-almost-every embedding $f$, or $(G, f)$ is not rigid for $\pi$-almost-every embedding. We call $G$ *rigid* if the former holds. Further details concerning this definition may be found in [13, 14, 27]. Note that rigid graphs are necessarily connected, but that there exist connected graphs which are not rigid.

We turn now to the rigidity of infinite graphs. Let $G$ be a countably infinite graph with finite vertex-degrees. The graph $G$ is called *rigid* if every finite subgraph of $G$ is contained in some finite rigid subgraph of $G$.

Next we introduce probability. Let $\mathcal{L}$ be a lattice in $d$ dimensions, and consider bond percolation on $\mathcal{L}$ having density $p$. The case $\mathcal{L} = \mathbb{L}^d$ is of no interest in the context of rigidity, since the lattice $\mathbb{L}^d$ is not itself rigid. Let $R$ be the event that the origin belongs to some infinite rigid subgraph of $\mathcal{L}$ all of whose edges are open. The *rigidity probability* is defined by

$$\theta^{\text{rig}}(p) = \mu_p(R).$$

Since $R$ is an increasing event, $\theta^{\text{rig}}$ is a non-decreasing function, whence

$$\theta^{\text{rig}}(p) \begin{cases} = 0 & \text{if } p < p_{\text{c}}^{\text{rig}}(\mathcal{L}), \\ > 0 & \text{if } p > p_{\text{c}}^{\text{rig}}(\mathcal{L}), \end{cases}$$

where the *rigidity critical probability* $p_{\text{c}}^{\text{rig}}(\mathcal{L})$ is given by

$$p_{\text{c}}^{\text{rig}}(\mathcal{L}) = \sup\{p : \theta^{\text{rig}}(p) = 0\}.$$

The study of the rigidity of percolation clusters was initiated by Jacobs and Thorpe, see [30, 31].

Since rigid graphs are connected, we have that $\theta^{\text{rig}}(p) \leq \theta(p)$, implying that $p_{\text{c}}^{\text{rig}}(\mathcal{L}) \geq p_{\text{c}}(\mathcal{L})$. The following is basic.

**Theorem 3 ([16, 27]).** *Let $\mathcal{L}$ be a $d$-dimensional lattice, where $d \geq 2$.*

(i) *We have that $p_{\text{c}}(\mathcal{L}) < p_{\text{c}}^{\text{rig}}(\mathcal{L})$.*
(ii) *$p_{\text{c}}^{\text{rig}}(\mathcal{L}) < 1$ if and only if $\mathcal{L}$ is rigid.*

How many (maximal) infinite rigid components may exist in a lattice $\mathcal{L}$? Let $J$ be the number of such components. By the Kolmogorov zero–one law, for any given value of $p$, $J$ is $\mu_p$-almost-surely constant. It may be conjectured that $\mu_p(J = 1) = 1$ whenever $\mu_p(J \geq 1) > 0$. The mathematical study of rigidity percolation was initiated by Holroyd in [27], where it was shown amongst other things that, for the triangular lattice $\mathbb{T}$ in two dimensions, $\mu_p(J = 1) = 1$ for almost every $p \in (p_{\text{c}}^{\text{rig}}(\mathbb{T}), 1]$. The proof was a highly non-trivial development of the Burton–Keane method. The main extra difficulty lies in the non-local nature of the property of rigidity. See also [29].

Considerably more general results have been obtained since by Häggström. Holroyd's result for *almost every $p$* was extended in [23] to *for every $p$*, by using the two-dimensional uniqueness arguments of Keane and co-authors to be found in [11]. More recently, Häggström has found an adaptation of the Burton–Keane argument which (almost) settles the problem for general rigid lattices in $d \geq 2$ dimensions.

**Theorem 4 ([25]).** *Let $d \geq 2$ and let $\mathcal{L}$ be a rigid $d$-dimensional lattice. We have that $\mu_p(J = 1) = 1$ whenever $p > p_{\text{c}}^{\text{rig}}(\mathcal{L})$.*

There remains the lacuna of deciding what happens when $p = p_{\text{c}}^{\text{rig}}$, that is, of proving either that $\mu_{p_{\text{c}}^{\text{rig}}}(J = 0) = 1$ or that $\mu_{p_{\text{c}}^{\text{rig}}}(J = 1) = 1$.

## 4. Entanglement in percolation

In addition to connectivity and rigidity, there is the notion of 'entanglement'. The simplest example of a graph which is entangled but not connected comprises two disjoint circuits which cannot be separated without one of them being broken. Such entanglement is intrinsically a three-dimensional affair, and therefore we restrict ourselves here to subgraphs of $\mathbb{L}^3 = (\mathbb{Z}^3, \mathbb{E})$ viewed as graphs embedded in a natural way in $\mathbb{R}^3$.

We begin with some terminology. For $E \subseteq \mathbb{E}$, we denote by $[E]$ the union of all unit line-segments of $\mathbb{R}^3$ corresponding to edges in $E$. The term 'sphere' is used to mean a subset of $\mathbb{R}^3$ which is homeomorphic to the 2-sphere $\{x \in \mathbb{R}^3 : \|x\| = 1\}$. The complement of any sphere $S$ has two connected components; we refer to the bounded component as the *inside* of $S$, written $\text{ins}(S)$, and to the unbounded component as the *outside* of $S$, written $\text{out}(S)$.

There is a natural definition of the term 'entanglement' when applied to finite sets of edges of the lattice $\mathbb{L}^3$, namely the following. We call the finite edge-set $E$ *entangled* if, for any sphere $S$ not intersecting $[E]$, either $[E] \subseteq \text{ins}(S)$ or $[E] \subseteq \text{out}(S)$. Thus entanglement is a property of edge-sets rather than of graphs. However, with any edge-set $E$ we may associate the graph $G_E$ having edge-set $E$ together with all incident vertices. Graphs $G_E$ arising in this way have no isolated vertices. We call $G_E$ *entangled* if $E$ is entangled, and we note that $G_E$ is entangled whenever it is connected.

There are several possible ways of extending the notion of entanglement to infinite subgraphs of $\mathbb{L}^3$, and these ways are not equivalent. For the sake of being definite, we adopt here a definition similar to that used for rigidity. Let $E$ be an infinite subset of $\mathbb{E}$. We call $E$ *entangled* if, for any finite subset $F$ ($\subseteq E$), there exists a finite entangled subset $F'$ of $E$ such that $F \subseteq F'$. We call the infinite graph $G_E$, defined as above, *entangled* if $E$ is entangled, and we note that $G_E$ is entangled whenever it is connected. A further discussion of the notion of entanglement may be found in [20].

Turning to percolation, we declare each edge of $\mathbb{L}^3$ to be open with probability $p$. We say that the origin 0 lies in an infinite open entanglement if there exists an infinite entangled set $E$ of open edges at least one of which has 0 as an endvertex. We concentrate on the event

$$N = \{0 \text{ lies in an infinite open entanglement}\},$$

and the *entanglement probability*

$$\theta^{\text{ent}}(p) = \mu_p(N).$$

Since $N$ is an increasing event, $\theta^{\text{ent}}$ is a non-decreasing function, whence

$$\theta^{\text{ent}}(p) \begin{cases} = 0 & \text{if } p < p_c^{\text{ent}}, \\ > 0 & \text{if } p > p_c^{\text{ent}}, \end{cases}$$

where the *entanglement critical probability* $p_c^{\text{ent}}$ is given by

$$p_c^{\text{ent}} = \sup\{p : \theta^{\text{ent}}(p) = 0\}. \tag{6}$$

Since every connected graph is entangled, it is immediate that $\theta(p) \leq \theta^{\text{ent}}(p)$, whence $0 \leq p_c^{\text{ent}} \leq p_c$.

**Theorem 5 ([1, 28]).** *The following strict inequalities are valid*:

$$0 < p_c^{\text{ent}} < p_c. \tag{7}$$

Entanglements in percolation appear to have been studied first in [32], where it was proposed that
$$p_c - p_c^{\text{ent}} \simeq 1.8 \times 10^{-7},$$
implying the strict inequality $p_c^{\text{ent}} < p_c$. It is a curious fact that we have no rigorous insight into the numerical value of $p_c^{\text{ent}}$. For example, we are unable to decide on the basis of mathematics whether $p_c^{\text{ent}}$ is numerically close to either 0 or $p_c$.

Suppose that $p$ is such that $\theta^{\text{ent}}(p) > 0$. By the zero–one law, the number $K$ of (maximal) infinite entangled open edge-sets satisfies $\mu_p(K \geq 1) = 1$. The almost-sure uniqueness of the infinite entanglement has been explored in [20, 29], and the situation is similar to that for rigidity percolation. Häggström's proof of the following theorem uses a non-trivial application of the Burton–Keane method.

**Theorem 6 ([24]).** *We have that $\mu_p(K = 1) = 1$ whenever $p > p_c^{\text{ent}}$.*

As is the case with rigidity, there remains the open problem of proving either that $\mu_{p_c^{\text{ent}}}(K = 0) = 1$ or that $\mu_{p_c^{\text{ent}}}(K = 1) = 1$.

There are several other open problems concerning entangled graphs, of which we mention a combinatorial question. Let $n \geq 1$, and let $\mathcal{E}_n$ be the set of all subsets $E$ of $\mathbb{E}$ with cardinality $n$ such that:

(i) some member of $e$ is incident to the origin, and
(ii) $E$ is entangled.

Since every connected graph is entangled, $\mathcal{E}_n$ is at least as large as the family of all connected sets of $n$ edges touching the origin. Therefore, $|\mathcal{E}_n|$ grows at least exponentially in $n$. It may be conjectured that there exists $\kappa$ such that
$$|\mathcal{E}_n| \leq e^{\kappa n} \qquad \text{for all } n \geq 1.$$
The best inequality known currently is of the form $|\mathcal{E}_n| \leq \exp\{\kappa n \log n\}$. See [20].

## 5. The random-cluster model

The random-cluster model on a finite graph $G = (V, E)$ is a certain parametric family of probability measures $\phi_{p,q}$ indexed by two parameters $p \in [0, 1]$ and $q \in (0, \infty)$. When $q = 1$, the measure is product measure with density $p$; when $q = 2, 3, \ldots$, the corresponding random-cluster measures correspond to the Ising and $q$-state Potts models on $G$. The random-cluster model provides a geometrical setting for the correlation functions of the ferromagnetic Ising and Potts models, and it has proved extremely useful in studying these models. Recent accounts of the theory, and of its impact on Ising/Potts models, may be found in [18, 19].

The configuration space is the set $\Omega = \{0, 1\}^E$ of 0/1-vectors indexed by the edge-set $E$. The probability measure $\phi_{p,q}$ on $\Omega$ is given by
$$\phi_{p,q}(\omega) = \frac{1}{Z} \left\{ \prod_{e \in E} p^{\omega(e)}(1-p)^{1-\omega(e)} \right\} q^{k(\omega)}, \qquad \omega \in \Omega, \tag{8}$$
where $k(\omega)$ is the number of connected components (or 'open clusters') of the graph $G_\omega = (V, \eta(\omega))$.

When $G$ is finite, every $\phi_{p,q}$-probability is a smooth function of the parameters $p$ and $q$. The situation is more interesting when $G$ is infinite, since infinite graphs may display phase transitions. For simplicity, we restrict the present discussion to the graph $\mathbb{L}^d = (\mathbb{Z}^d, \mathbb{E})$ where $d \geq 2$. We introduce next the concept of boundary conditions.

Let $\Lambda$ be a finite box, and write $\mathbb{E}_\Lambda$ for the set of edges joining pairs of members of $\Lambda$. We write $\mathcal{T}_\Lambda$ for the $\sigma$-field generated by the states of edges in $\mathbb{E} \setminus \mathbb{E}_\Lambda$. For $\xi \in \Omega$, we write $\Omega^\xi_\Lambda$ for the (finite) subset of $\Omega$ containing all configurations $\omega$ satisfying $\omega(e) = \xi(e)$ for $e \in \mathbb{E}^d \setminus \mathbb{E}_\Lambda$; these are the configurations which 'agree with $\xi$ off $\Lambda$'. Let $\xi \in \Omega$, and write $\phi^\xi_{\Lambda,p,q}$ for the random-cluster measure on the finite graph $\Lambda$ 'with boundary condition $\xi$'. That is to say, $\phi_{\Lambda,p,q}$ is given as in (8) subject to $\omega \in \Omega^\xi_\Lambda$, and with $k(\omega)$ replaced by the number of open clusters of $\mathbb{L}^d$ that intersect $\Lambda$.

A probability measure $\phi$ on $(\Omega, \mathcal{F})$ is called a *random-cluster measure* with parameters $p$ and $q$ if for all $A \in \mathcal{F}$ and all finite boxes $\Lambda$, $\quad \phi(A \mid \mathcal{T}_\Lambda)(\xi) = \phi^\xi_{\Lambda,p,q}(A) \quad$ for $\phi$-a.e. $\xi$.

The set of such measures is denoted $\mathcal{R}_{p,q}$. The reader is referred to [18, 19] for accounts of the existence and properties of random-cluster measures.

One may construct infinite-volume measures by taking weak limits. A probability measure $\phi$ on $(\Omega, \mathcal{F})$ is called a *limit random-cluster measure* with parameters $p$ and $q$ if there exist $\xi \in \Omega$ and a sequence $\mathbf{\Lambda} = (\Lambda_n : n \geq 1)$ of boxes satisfying $\Lambda_n \to \mathbb{Z}^d$ as $n \to \infty$ such that

$$\phi^\xi_{\Lambda_n,p,q} \Rightarrow \phi \qquad \text{as } n \to \infty.$$

The two 'extremal' boundary conditions are the configurations 'all 0' and 'all 1', denoted by 0 and 1, respectively. It is a standard application of positive association that the weak limits

$$\phi^b_{p,q} = \lim_{\Lambda \uparrow \mathbb{Z}^d} \phi^b_{\Lambda,p,q}$$

exist for $b = 0, 1$ and $q \geq 1$. It is an important fact that these limits belong to $\mathcal{R}_{p,q}$.

**Theorem 7 ([15]).** *Let $p \in [0,1]$ and $q \in [1, \infty)$. The limit random-cluster measures $\phi^b_{p,q}$, $b = 0, 1$, belong to $\mathcal{R}_{p,q}$.*

The proof hinges on the following fact. Let $\phi$ be a limit random-cluster measure with parameters $p$, $q$ such that the number $I$ of infinite open clusters satisfies

$$\phi(I \in \{0, 1\}) = 1. \tag{9}$$

It may then be deduced that $\phi \in \mathcal{R}_{p,q}$. The uniqueness theorem, Theorem 2, is used to establish (9) for the measures $\phi = \phi^b_{p,q}$, $b = 0, 1$.

Let $q \geq 1$. The random-cluster model has a phase transition defined as follows. For $b = 0, 1$, let $\theta^b(p, q) = \phi^b_{p,q}(0 \leftrightarrow \infty)$, and define the critical points

$$p^b_c(q) = \sup\{p : \theta^b(p, q) = 0\}.$$

It is standard that $\phi^0_{p,q} = \phi^1_{p,q}$ for almost every $p$. It follows that $p^0_c(q) = p^1_c(q)$, and we write $p_c(q)$ for the common critical value. It is known that $\phi^0_{p,q} = \phi^1_{p,q}$ when $p < p_c(q)$, and it is an important open problem to prove that

$$\phi^0_{p,q} = \phi^1_{p,q} \qquad \text{if } p > p_c(q).$$

See [18, 19] for further discussion of the uniqueness of random-cluster measures.

## 6. Random walks in random labyrinths

Suppose that a ball is propelled through a random environment of obstacles, off which it rebounds with perfect reflection. We ask for information about the trajectory of the ball. This classical problem is often termed the 'Lorentz problem', and it has received much attention in both the mathematics and physics literature. If the obstacles are distributed at random in $\mathbb{R}^d$ then, conditional on their placements, the motion of the ball is deterministic. It is a significant problem of probability theory to develop a rigorous analysis of such a situation.

Two natural questions spring to mind.

(i) *Non-localization.* What is the probability that the trajectory of the ball is unbounded?
(ii) *Diffusivity.* Suppose the trajectory is unbounded with a strictly positive probability. Conditional on this event, is there a central limit theorem for the ball's position after a large time $t$.

These questions seem to be difficult, especially when the obstacles are distributed aperiodically. The problem is much easier when the environment of obstacles is 'lubricated' by a positive density of points at which the ball behaves as a random walk.

We consider a lattice model of the following type. The obstacles are distributed around the vertex-set of the $d$-dimensional hypercubic lattice $\mathbb{L}^d$, and they are designed in such a way that the ball traverses the edges of the lattice. Some of the associated mathematics has been surveyed in [4, 17], to which the reader is referred for an account of the literature. The main result of [4] is that, subject to certain conditions on the density of obstacles, the ball's trajectory satisfies a functional central limit theorem. The Burton–Keane method plays a crucial role in the proof.

We make this more concrete as follows. Our model involves a random environment of reflecting bodies distributed around the vertices of $\mathbb{L}^d$. Each vertex is designated either a 'reflector' (of a randomly chosen type) or a 'random walk point'. The interpretation of the term 'random walk point' is as follows: when the ball hits such a point, it departs in a direction chosen randomly from the $2d$ available directions, this direction being chosen independently of everything else.

The defining properties of a reflector $\rho$ are that:

(i) to each incoming direction $u$ there is assigned a unique outgoing direction $\rho(u)$, and
(ii) the ball will retrace its path if its direction is reversed.

Let $I = \{e_1, e_2, \ldots, e_d\}$ be the set of positive unit vectors of $\mathbb{Z}^d$, and let $I^\pm = \{\alpha e_j : \alpha = \pm,\ 1 \leq j \leq d\}$. A *reflector* is defined to be a map $\rho : I^\pm \to I^\pm$ with the property that
$$\rho(-\rho(u)) = -u \qquad \text{for all } u \in I^\pm.$$

We write $\mathcal{R}$ for the set of all reflectors. One particular reflector is special, namely the identity map satisfying $\rho(u) = u$ for all $u \in I^\pm$; we call this the *crossing*, and we denote it by $+$. Crossings do not deflect the ball.

The following random environment will be termed a *random labyrinth*. Let $p_{\text{rw}}$ and $p_+$ be non-negative reals such that $p_{\text{rw}} + p_+ \leq 1$, and let $\pi$ be a probability mass function on the set $\mathcal{R} \setminus \{+\}$ of 'non-trivial' reflectors (that is, $\pi(\rho) \geq 0$ for $\rho \in \mathcal{R}\setminus\{+\}$ and $\sum_{\rho \in \mathcal{R}\setminus\{+\}} \pi(\rho) = 1$). Let $Z = (Z_x : x \in \mathbb{Z}^d)$ be a family of independent random variables, taking values in $\mathcal{R} \cup \{\varnothing\}$, with probabilities

$$\mathbb{P}(Z_x = \beta) = \begin{cases} p_{\mathrm{rw}} & \text{if } \beta = \varnothing, \\ p_+ & \text{if } \beta = +, \\ (1 - p_{\mathrm{rw}} - p_+)\pi(\rho) & \text{if } \beta = \rho \in \mathcal{R}\setminus\{+\}. \end{cases}$$

A vertex $x$ is called a *crossing* if $Z_x = +$, and a *random walk (rw) point* if $Z_x = \varnothing$.

We now study admissible paths in the labyrinth $Z$. Consider a path in $\mathbb{L}^d$ which visits (in order) the vertices $x_0, x_1, \ldots, x_n$; we allow the path to revisit a given vertex more than once, and to traverse a given edge more than once. This path is called *admissible* if it conforms to the reflectors that it meets, which is to say that

$$x_{j+1} - x_j = Z_{x_j}(x_j - x_{j-1}) \qquad \text{for all } j \text{ such that } Z_{x_j} \neq \varnothing.$$

Remarkably little is known about random labyrinths when $p_{\mathrm{rw}} = 0$. One notorious open problem concerns the existence (or not) of infinite admissible paths in $\mathbb{L}^2$ when $p_{\mathrm{rw}} = 0$. The problem is substantially easier when $p_{\mathrm{rw}} > 0$, and we assume this henceforth. We explain next how the labyrinth $Z$ generates a random walk therein. Let $x$ be a rw point. A walker starts at $x$, and flips a fair $2d$-sided coin in order to determine the direction of its first step. Henceforth, it is required to traverse admissible paths only, and it flips the coin to determine its exit direction from any rw point encountered. We write $P_x^Z$ for the law of the random walk in the labyrinth $Z$, starting from a rw point $x$.

There is a natural equivalence relation on the set $R$ of rw points of $\mathbb{Z}^d$, namely $x \leftrightarrow y$ if there exists an admissible path with endpoints $x$ and $y$. Let $C_x$ be the equivalence class containing the rw point $x$. We may follow the progress of a random walk starting at $x$ by writing down (in order) the rw points which it visits, say $X_0 (= x), X_1, X_2, \ldots$. Given the labyrinth $Z$, the sequence $X = (X_n : n \geq 0)$ is an irreducible Markov chain on the countable state space $C_x$. Furthermore, this chain is reversible with respect to the measure $\mu$ given by $\mu(y) = 1$ for $y \in C_x$. We say that $x$ is *$Z$-localized* if $|C_x| < \infty$, and *$Z$-non-localized* otherwise. We call $Z$ *localized* if all rw points are $Z$-localized, and we call $Z$ *non-localized* otherwise. By a zero–one law, we have that $\mathbb{P}(Z \text{ is localized})$ equals either 0 or 1.

Suppose that the origin 0 is a rw point. As before, we consider the sequence $X_0 (= 0), X_1, X_2, \ldots$ of rw points visited in sequence by a random walk in $Z$ beginning at the origin 0. For $\epsilon > 0$, we let

$$X^\epsilon(t) = \epsilon X_{\lfloor \epsilon^{-2} t \rfloor} \qquad \text{for } t \geq 0,$$

and we are interested in the behaviour of the process $X^\epsilon(\cdot)$ in the limit as $\epsilon \downarrow 0$. We study $X^\epsilon$ under the probability measure $\mathbb{P}_0$, defined as the measure $\mathbb{P}$ conditional on the event $\{0 \text{ is a rw point, and } |C_0| = \infty\}$.

We write $p_{\mathrm{c}}^{\mathrm{site}}$ for the critical probability of site percolation on $\mathbb{Z}^d$.

**Theorem 8 ([4]).** *Let $d \geq 2$ and $p_{\mathrm{rw}} > 0$. There exists a strictly positive constant $A = A(p_{\mathrm{rw}}, d)$ such that the following holds whenever either $1 - p_{\mathrm{rw}} - p_+ < A$ or $p_{\mathrm{rw}} > p_{\mathrm{c}}^{\mathrm{site}}$:*

(i) $\mathbb{P}(0 \text{ is a rw point, and } |C_0| = \infty) > 0$, *and*
(ii) *as $\epsilon \downarrow 0$, the re-scaled process $X^\epsilon(\cdot)$ converges $\mathbb{P}_0$-dp to $\sqrt{\delta}W$, where $W$ is a standard Brownian motion in $\mathbb{R}^d$ and $\delta$ is a strictly positive constant.*

With $E$ denoting expectation, the convergence '$\mathbb{P}_0$-dp' is to be interpreted as

$$P_0^Z(f(X^\epsilon)) \to E(f(W)) \qquad \text{in } \mathbb{P}_0\text{-probability},$$

for all bounded continuous functions $f$ on the Skorohod path-space $D([0,\infty), \mathbb{R}^d)$.

The proof of Theorem 8 utilizes the Kipnis–Varadhan central limit theorem, [34], together with its application to percolation, see [7, 8]. A key step in the proof is to show that, under the conditions of the theorem, there exists a *unique* infinite equivalence class, and this is where the Burton–Keane method is key.

## 7. Non-uniqueness for non-amenable graphs

Let $G = (V,E)$ be an infinite, connected graph with finite vertex-degrees. We call $G$ *amenable* if its 'isoperimetric constant'

$$\chi(G) = \inf\left\{\frac{|\partial W|}{|W|} : W \subseteq V,\ 0 < |W| < \infty\right\} \tag{10}$$

satisfies

$$\chi(G) = 0,$$

where the infimum in (10) is over all non-empty finite subsets $W$ of $V$. The graph is called *non-amenable* if $\chi(G) > 0$.

We have so far concentrated on situations where infinite clusters are (almost surely) unique, as is commonly the case for an amenable graph. The situation is quite different when the graph is non-amenable, and a systematic study of percolation on such graphs was proposed in [3]. The best known example is bond percolation on the infinite binary tree, for which there exist infinitely many infinite clusters whenever the edge-density $p$ satisfies $\frac{1}{2} < p < 1$. Let $G = (V,E)$ be an infinite graph and let $p \in (0,1)$. For $\omega \in \Omega = \{0,1\}^E$, let $I = I(\omega)$ be the number of infinite clusters of $\omega$. It has been known since [22] that there exist graphs having three non-trivial phases characterized respectively by $I = 0$, $I = 1$, $I = \infty$. One of the most interesting results in this area is the existence of a critical point for the event $\{I = 1\}$. This is striking because the event $\{I = 1\}$ is not increasing. Prior to stating this theorem, we introduce some jargon.

The infinite connected graph $G = (V,E)$ is called *transitive* if, for all $x,y \in V$, there exists an automorphism $\tau$ of $G$ such that $y = \tau(x)$. The graph $G$ is called *quasi-transitive* if there exists a finite set $V_0$ of vertices such that, for all $y \in V$, there exists $x \in V_0$ and an automorphism $\tau$ such that $y = \tau(x)$.

The following result was obtained by Häggström and Peres under a further condition, subseqently lifted by Schonmann.

**Theorem 9 ([26, 35]).** *Let $G$ be an infinite, connected, quasi-transitive graph. There exist $p_c, p_u \in [0,1]$ satisfying $0 \leq p_c \leq p_u \leq 1$ such that:*

$$\mu_p(I = 0) = 1 \quad \text{if} \quad 0 \leq p < p_c, \tag{11}$$
$$\mu_p(I = \infty) = 1 \quad \text{if} \quad p_c < p < p_u, \tag{12}$$
$$\mu_p(I = 1) = 1 \quad \text{if} \quad p_u < p \leq 1. \tag{13}$$

Here are some examples.

1. For an amenable graph, we have by the Burton–Keane argument that $p_c = p_u$.
2. For the binary tree, we have $p_c = \frac{1}{2}$ and $p_u = 1$.
3. For the direct product of the binary tree and a line, we have that $0 < p_c < p_u < 1$, see [22].